\documentclass{ws-m3as}

\usepackage{graphicx}
\usepackage{psfrag}
\usepackage[latin1]{inputenc}

\newcommand{\real}{\mathbb{R}}

\begin{document}


\markboth{J. Jansson, C. Johnson, and A. Logg}{Computational Modeling of Dynamical Systems}

\title{COMPUTATIONAL MODELING OF DYNAMICAL SYSTEMS}

\author{JOHAN JANSSON}
\address{Department of Computational Mathematics,
         Chalmers University of Technology, SE--412 96 G\"{o}teborg, Sweden,
	 \emph{email}: johanjan{\@@}math.chalmers.se.}

\author{CLAES JOHNSON}
\address{Department of Computational Mathematics,
         Chalmers University of Technology, SE--412 96 G\"{o}teborg, Sweden,
	 \emph{email}: claes{\@@}math.chalmers.se.}

\author{ANDERS LOGG}
\address{Department of Computational Mathematics,
         Chalmers University of Technology, SE--412 96 G\"{o}teborg, Sweden,
	 \emph{email}: logg{\@@}math.chalmers.se.}

\maketitle

\begin{abstract}
  In this short note, we discuss the basic approach to computational
  modeling of dynamical systems. If a dynamical system contains
  multiple time scales, ranging from very fast to slow, computational
  solution of the dynamical system can be very costly.
  By resolving the fast time scales in a short time simulation, a
  model for the effect of the small time scale variation on large time
  scales can be determined, making solution possible on a long time interval.
  This process of computational modeling can be completely automated.
  Two examples are presented, including a simple model problem
  oscillating at a time scale of $10^{-9}$ computed over the
  time interval $[0,100]$, and a lattice
  consisting of large and small point masses.
\end{abstract}


\section{Introduction}

We consider a dynamical system of the form
\begin{equation}
  \label{eq:u'=f}
  \begin{array}{rcl}
    \dot{u}(t) &=& f(u(t),t), \quad t \in (0,T], \\
    u(0) &= u_0,
  \end{array}
\end{equation}
where $u : [0,T] \rightarrow \real^N$ is the solution to be computed,
$u_0 \in \real^N$ a given initial value, $T>0$ a given final time,
and $f : \real^N \times (0,T] \rightarrow \real^N$ a given function
that is Lipschitz-continuous in $u$ and bounded. We consider a
situation where the exact solution $u$ varies on different time
scales, ranging from very fast to slow. Typical examples include
meteorological models for weather prediction, with fast time scales on
the range of seconds and slow time scales on the range of years,
protein folding represented by a molecular dynamics model of
the form (\ref{eq:u'=f}), with fast
time scales on the range of femtoseconds and slow time scales on the
range of microseconds, or turbulent flow with a wide range of time scales.

To make computation feasible in a situation where computational
resolution of the fast time scales would be prohibitive because of the
small time steps required, the given model (\ref{eq:u'=f}) containing
the fast time scales needs to be replaced with a \emph{reduced model}
for the variation of the solution $u$ of (\ref{eq:u'=f}) on resolvable time scales.
As discussed below, the key step is to correctly model the effect
of the variation at the fast time scales on the variation on slow time scales.

The problem of model reduction is very general and various approaches
have been taken\cite{RuhSko98,Kre91}. We
present below a new approach to model reduction, based on resolving
the fast time scales in a short time simulation and determining a
model for the effect of the small time scale variation on large time
scales. This process of computational modeling can be completely
\emph{automated} and the validity of the reduced model can be
evaluated a posteriori.

\section{A simple model problem}

\begin{figure}[htbp]
  \begin{center}
    \psfrag{k1}{$\kappa = 1$}
    \psfrag{k1}{$\kappa = 1$}
    \psfrag{k2}{$\kappa \gg 1$}
    \psfrag{u1}{$u_1$}
    \psfrag{u2}{$u_2$}
    \psfrag{M1}{$M_1$}
    \psfrag{M2}{$M_2$}
    \includegraphics[width=10cm]{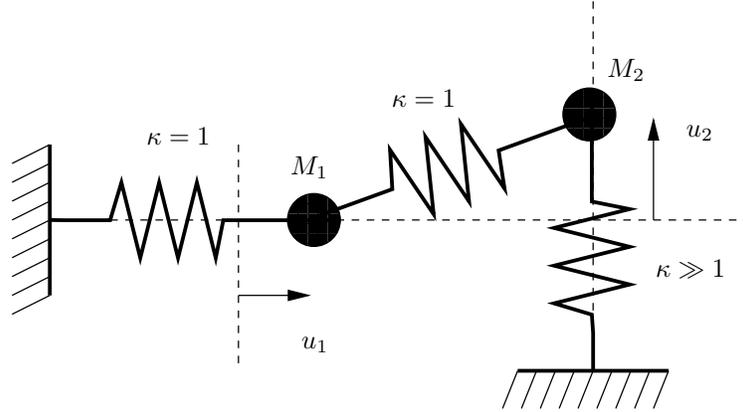}
    \caption{A simple mechanical system with large time scale $\sim 1$ and
      small time scale $\sim 1/\sqrt{\kappa}$.}
    \label{fig:model}
  \end{center}
\end{figure}

We consider a simple example illustrating the basic aspects:
Find $u=(u_1,u_2):[0,T] \rightarrow \real^2$, such that
\begin{equation}
  \label{eq:model}
  \begin{array}{rcl}
    \ddot{u}_1 + u_1 - u_2^2/2 &=& 0 \quad \mbox{on } (0,T],\\
      \ddot{u}_2 + \kappa u_2 &=& 0 \quad \mbox{on } (0,T],\\
	u(0) = (0,1) &&\quad \dot{u}(0) = (0,0),
  \end{array}
\end{equation}
which models a moving unit point mass $M_1$ connected
through a soft spring to another unit point mass $M_2$, with
$M_2$ moving along a line perpendicular to the line of motion of $M_1$,
see Figure \ref{fig:model}. The second point mass $M_2$ is connected
to a fixed support through a very stiff spring with
spring constant $\kappa = 10^{18}$ and oscillates rapidly on a time
scale of size $1/\sqrt{\kappa} = 10^{-9}$.
The oscillation of $M_2$ creates a force $\sim u_2^2$ on $M_1$
proportional to the elongation of the spring connecting $M_2$ to $M_1$
(neglecting terms of order $u_2^4$).

The short time scale of size $10^{-9}$ requires
time steps of size $\sim 10^{-10}$ for full resolution. With $T = 100$,
this means a total of $\sim 10^{12}$ time steps for solution of (\ref{eq:model}).
However, by replacing (\ref{eq:model}) with a reduced model where the
fast time scale has been removed, it is possible to compute the (averaged)
solution of (\ref{eq:model}) with time steps of size $\sim 0.1$
and consequently only a total of $10^3$ time steps.

\section{Taking averages to obtain the reduced model}

Having realized that point-wise resolution of the fast time scales of
the exact solution $u$ of (\ref{eq:u'=f}) may sometimes be computationally very
expensive or even impossible, we seek instead to compute a time average
$\bar{u}$ of $u$, defined by
\begin{equation}
  \label{eq:average}
  \bar{u}(t)= \frac{1}{\tau}\int_{-\tau/2}^{\tau/2}u(t+s)\, ds, \quad
  t \in [\tau/2,T - \tau/2],
\end{equation}
where $\tau > 0$ is the size of the average.
The average $\bar{u}$ can be extended to $[0,T]$ in various ways.
We consider here a constant extension, i.e.,
we let $\bar{u}(t) = \bar{u}(\tau/2)$ for $t \in [0,\tau/2)$,
and let $\bar{u}(t) = \bar{u}(T-\tau/2)$ for $t \in (T - \tau/2,T]$.

We now seek a dynamical system satisfied by the average $\bar{u}$
by taking the average of (\ref{eq:u'=f}). We obtain
\begin{displaymath}
  \dot{\bar{u}}(t) = \bar{\dot{u}}(t) = \overline{f(u,\cdot)}(t)
  = f(\bar{u}(t),t) + (\overline{f(u,\cdot)}(t) - f(\bar{u}(t),t)),
\end{displaymath}
or
\begin{equation}
  \label{eq:u'=f,average}
  \dot{\bar{u}}(t) = f(\bar{u}(t),t) + \bar{g}(u,t),
\end{equation}
where the \emph{variance} $\bar{g}(u,t) = \overline{f(u,\cdot)}(t) - f(\bar{u}(t),t)$
accounts for the effect of small scales on time scales larger than $\tau$.
(Note that we may extend (\ref{eq:u'=f,average}) to $(0,T]$ by
defining $\bar{g}(u,t) = - f(\bar{u}(t),t)$ on $(0,\tau/2] \cup (T-\tau/2,T]$.)

We now seek to model the variance $\bar{g}(u,t)$ in the form
$\bar{g}(u,t)\approx \tilde{g}(\bar{u}(t),t)$ and replace (\ref{eq:u'=f,average}) and
thus (\ref{eq:u'=f}) by
\begin{equation}
  \label{eq:reduced}
  \begin{array}{rcl}
    \dot{\tilde{u}}(t) &=& f(\tilde{u}(t),t) + \tilde{g}(\tilde{u}(t),t),
    \quad t \in (0,T], \\
      \tilde{u}(0) &=& \bar{u}_0,
  \end{array}
\end{equation}
where $\bar{u}_0 = \bar{u}(0) = \bar{u}(\tau/2)$.
We refer to this system as the \emph{reduced model} with \emph{subgrid
model} $\tilde{g}$ corresponding to (\ref{eq:u'=f}).

To summarize, if the solution $u$ of the full dynamical system (\ref{eq:u'=f})
is computationally unresolvable, we aim at computing the
average $\bar{u}$ of $u$. However, since the
variance $\bar{g}$ in the averaged dynamical system (\ref{eq:u'=f,average}) is
unknown, we need to solve the reduced model (\ref{eq:reduced})
for $\tilde{u} \approx \bar{u}$ with an approximate subgrid
model $\tilde{g} \approx \bar{g}$. Solving the
reduced model (\ref{eq:reduced}) using e.g. a Galerkin finite element
method, we obtain an approximate solution $U \approx \tilde{u} \approx \bar{u}$.
Note that we may not expect $U$ to be close to $u$ point-wise in time,
while we hope that $U$ is close to $\bar{u}$ point-wise.

\section{Modeling the variance}

There are two basic approaches to the modeling of the variance $\bar{g}(u,t)$ in
the form $\tilde{g}(\tilde{u}(t),t)$; (i) scale-extrapolation or (ii) local resolution.
In (i), a sequence of solutions is computed with increasingly fine
resolution, but without resolving the fastest time scales. A model for
the effects of the fast unresolvable scales is then determined by
extrapolation from the sequence of computed solutions\cite{Hof02}.
In (ii), the approach followed below, the solution $u$ is computed
accurately over a short time period, resolving the fastest time
scales. The reduced model is then obtained by computing the variance
\begin{equation}
  \bar{g}(u,t) = \overline{f(u,\cdot)}(t) - f(\bar{u}(t),t)
\end{equation}
and then determining $\tilde{g}$ for the remainder of the time interval
such that $\tilde{g}(\tilde{u}(t),t) \approx \bar{g}(u,t)$.

For the simple model problem (\ref{eq:model}), which we can write in the
form (\ref{eq:u'=f}) by introducing the two new variables $u_3 = \dot{u}_1$
and $u_4 = \dot{u}_2$ with
\begin{displaymath}
  f(u,\cdot) = (u_3, u_4, -u_1 + u_2^2/2, -\kappa u_2),
\end{displaymath}
we note
that $\bar{u}_2 \approx 0$ (for $\sqrt{\kappa} \tau$ large)
while $\overline{u_2^2} \approx 1/2$.
By the linearity of $f_1$, $f_2$, and $f_4$, the (approximate) reduced
model takes the form
\begin{equation}
  \label{eq:model,reduced}
  \begin{array}{rcl}
    \ddot{\tilde{u}}_1 + \tilde{u}_1 - 1/4 &=& 0 \quad \mbox{on } (0,T],\\
    \ddot{\tilde{u}}_2 + \kappa \tilde{u}_2 &=& 0 \quad \mbox{on } (0,T],\\
      \tilde{u}(0) = (0,0), &&\quad \dot{\tilde{u}}(0) = (0,0),
  \end{array}
\end{equation}
with solution $\tilde{u}(t) = (\frac{1}{4}(1 - \cos t),0)$.

In general, the reduced model is constructed with subgrid model $\tilde{g}$
varying on resolvable time scales. In the simplest case, it is enough to
model $\tilde{g}$ with a constant and repeatedly checking the validity
of the model by comparing the reduced model (\ref{eq:reduced})
with the full model (\ref{eq:u'=f}) in a short time simulation. Another
possibility is to use a piecewise polynomial representation for the
subgrid model $\tilde{g}$.

\section{Solving the reduced system}

Although the presence of small scales has been decreased in the
reduced system (\ref{eq:reduced}), the small scale variation may still
be present. This is not evident in the reduced system (\ref{eq:model,reduced})
for the simple model problem (\ref{eq:model}), where we
made the approximation $\tilde{u}_2(0) = 0$. In practice, however, we
compute $\tilde{u}_2(0) = \frac{1}{\tau}\int_0^\tau u_2(t) \, dt =
\frac{1}{\tau}\int_0^\tau \cos(\sqrt{\kappa} t) \, dt \sim
1/(\sqrt{\kappa}\tau)$ and so $\tilde{u}_2$ oscillates at the fast
time scale $1/\sqrt{\kappa}$ with amplitude $1/(\sqrt{\kappa}\tau)$.

To remove these oscillations, the reduced system needs to be stabilized by
introducing damping of high frequencies. Following the general
approach\cite{HofJoh04b}, a least squares stabilization is
added in the Galerkin formulation of the reduced system (\ref{eq:reduced})
in the form of a modified test function. As a result, damping is
introduced for high frequencies without affecting low frequencies.

Alternatively, components such as $u_2$ in (\ref{eq:model,reduced})
may be \emph{inactivated}, corresponding to a subgrid model of
the form $\tilde{g}_2(\tilde{u},\cdot) = -f_2(\tilde{u},\cdot)$.
We take this simple approach for the example problems presented below.

\section{Error analysis}

The validity of a proposed subgrid model may be checked a posteriori.
To analyze the modeling error introduced by approximating the
variance $\bar{g}$ with the subgrid model $\tilde{g}$,
we introduce the \emph{dual problem}
\begin{equation}
  \label{eq:dual}
  \begin{array}{rcl}
    - \dot{\phi}(t) &=& J(\bar{u},U,t)^{\top} \phi(t), \quad t \in [0,T), \\
      \phi(T) &=& \psi,
  \end{array}
\end{equation}
where $J$ denotes the Jacobian of the right-hand side of the
dynamical system (\ref{eq:u'=f}) evaluated at a mean value of the
average $\bar{u}$ and the computed numerical (finite element)
solution $U \approx \tilde{u}$ of the reduced system (\ref{eq:reduced}),
\begin{equation}
  J(\bar{u},U,t) = \int_0^1 \frac{\partial f}{\partial u}
  (s \bar{u}(t) + (1-s) U(t),t) \, ds,
\end{equation}
and where $\psi$ is initial data for the backward dual problem.

To estimate the error $\bar{e} = U - \bar{u}$ at final time,
we note that $\bar{e}(0) = 0$
and $\dot{\phi} + J(\bar{u},U,\cdot)^{\top} \phi = 0$, and write
\begin{displaymath}
  \begin{array}{rcl}
    (\bar{e}(T),\psi)
    &=& (\bar{e}(T),\psi) -
    \int_0^T (\dot{\phi} + J(\bar{u},U,\cdot)^{\top} \phi,
    \bar{e}) \, dt \\
    &=& \int_0^T (\phi, \dot{\bar{e}} - J\bar{e}) \, dt
    = \int_0^T (\phi, \dot{U} - \dot{\bar{u}} -
    f(U,\cdot) + f(\bar{u},\cdot)) \, dt \\
    &=& \int_0^T (\phi, \dot{U} - f(U,\cdot) - \tilde{g}(U,\cdot)) \, dt +
    \int_0^T (\phi, \tilde{g}(U,\cdot) - \bar{g}(u,\cdot)) \, dt \\
    &=& \int_0^T (\phi, \tilde{R}(U,\cdot)) \, dt +
    \int_0^T (\phi, \tilde{g}(U,\cdot) - \bar{g}(u,\cdot)) \, dt.
  \end{array}
\end{displaymath}
The first term, $\int_0^T (\phi, \tilde{R}(U,\cdot)) \, dt$, in this
\emph{error representation} corresponds to
the \emph{discretization error} $U - \tilde{u}$
for the numerical solution of (\ref{eq:reduced}). If a Galerkin finite
element method is used\cite{EriEst95,EriEst96},
the \emph{Galerkin orthogonality} expressing the orthogonality of the
residual $\tilde{R}(U,\cdot) = \dot{U} - f(U,\cdot) - \tilde{g}(U,\cdot)$
to a space of test functions can be used to subtract a
test space interpolant $\pi \phi$ of the dual solution $\phi$.
In the simplest case of the $\mathrm{cG}(1)$ method for a partition
of the interval $(0,T]$ into $M$ subintervals $I_j = (t_{j-1},t_j]$, each of
length $k_j = t_j - t_{j-1}$, we subtract a
piecewise constant interpolant to obtain
\begin{displaymath}
  \begin{array}{rcl}
    \int_0^T (\phi, \tilde{R}(U,\cdot)) \, dt
    &=& \int_0^T (\phi - \pi \phi, \tilde{R}(U,\cdot)) \, dt
    \leq \sum_{j=1}^M k_j \max_{I_j} \|\tilde{R}(U,\cdot)\|_{l_2}
    \int_{I_j} \|\dot{\phi}\|_{l_2} \, dt \\
    &\leq& S^{[1]}(T) \max_{[0,T]} \|k\tilde{R}(U,\cdot)\|_{l_2},
  \end{array}
\end{displaymath}
where the \emph{stability factor} $S^{[1]}(T) = \int_0^T \|\dot{\phi}\|_{l_2} \, dt$
measures the sensitivity to discretization errors for the given output
quantity $(\bar{e}(T),\psi)$.

The second term, $\int_0^T (\phi, \tilde{g}(U,\cdot) - \bar{g}(u,\cdot)) \, dt$,
in the error representation corresponds
to the \emph{modeling error} $\tilde{u} - \bar{u}$. The sensitivity to
modeling errors is measured by the stability
factor $S^{[0]}(T) = \int_0^T \|\phi\|_{l_2} \, dt$.
We notice in particular that if the stability factor $S^{[0]}(T)$ is of moderate
size, a reduced model of the form (\ref{eq:reduced})
for $\tilde{u} \approx \bar{u}$ may be constructed.

We thus obtain
the error estimate
\begin{equation}
  \vert (\bar{e}(T),\psi) \vert
  \leq S^{[1]}(T) \max_{[0,T]} \|k\tilde{R}(U,\cdot)\|_{l_2} +
  S^{[0]}(T) \max_{[0,T]} \|\tilde{g}(U,\cdot) - \bar{g}(u,\cdot)\|_{l_2},
\end{equation}
including both discretization and modeling errors. The initial
data $\psi$ for the dual problem (\ref{eq:dual}) is chosen to reflect the desired
output quantity, e.g. $\psi = (1,0,\ldots,0)$ to measure the error
in the first component of $U$.

To estimate the modeling error, we need to estimate the
quantity $\tilde{g} - \bar{g}$.
This estimate is obtained by repeatedly solving the full dynamical
system (\ref{eq:u'=f}) at a number of control points and comparing the
subgrid model $\tilde{g}$ with the computed variance $\bar{g}$. As
initial data for the full system at a control point, we
take the computed solution $U \approx \bar{u}$ at the control point
and add a perturbation of appropriate size, with the size of the
perturbation chosen to reflect the initial oscillation at the fastest
time scale.

\section{Numerical results}

We present numerical results for two model problems, including the
simple model problem (\ref{eq:model}), computed with DOLFIN\cite{logg:www:01} version
0.4.10. With the option \emph{automatic modeling}
set, DOLFIN automatically creates the reduced model (\ref{eq:reduced})
for a given dynamical system of the form (\ref{eq:u'=f}) by resolving the full
system in a short time simulation and then determining
a constant subgrid model $\bar{g}$. Components with constant average,
such as $u_2$ in (\ref{eq:model}), are automatically marked as
inactive and are kept constant throughout the simulation. The
automatic modeling implemented in DOLFIN is rudimentary and many
improvements are possible, but it represents a first attempt at the
automation of modeling, following the recently presented\cite{logg:thesis:03}
directions for the \emph{automation of computational mathematical modeling}.

\subsection{The simple model problem}

The solution for the two components of the simple model problem (\ref{eq:model})
is shown in Figure \ref{fig:model,solution} for $\kappa = 10^{18}$
and $\tau = 10^{-7}$. The value of the
subgrid model $\bar{g}_1$ is automatically determined to $0.2495 \approx 1/4$.

\begin{figure}[htbp]
  \begin{center}
    \psfrag{t}{$t$}
    \psfrag{u1}{$u_1(t)$}
    \psfrag{u2}{$u_2(t)$}
    \includegraphics[width=10cm]{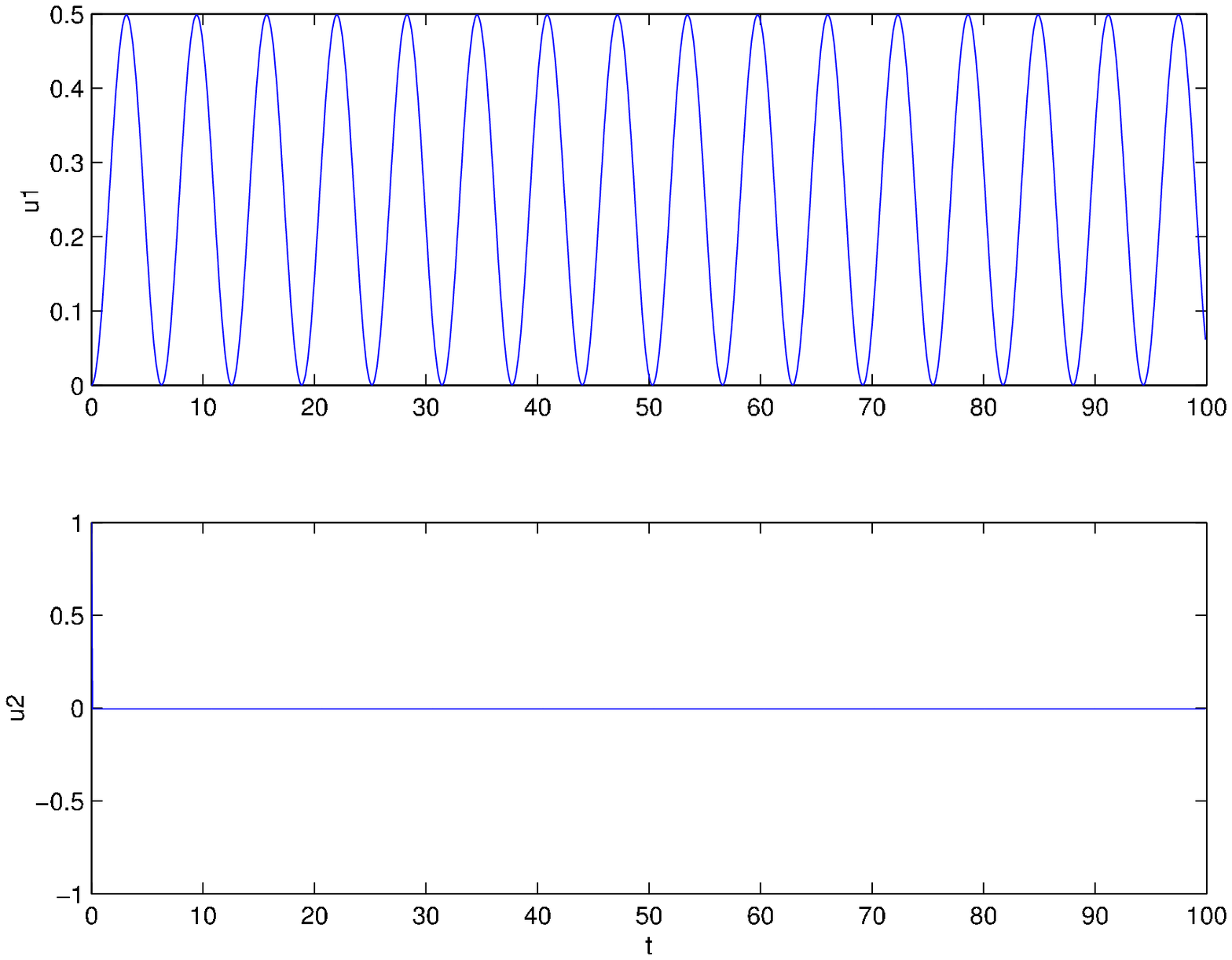} \\
    \includegraphics[width=10cm]{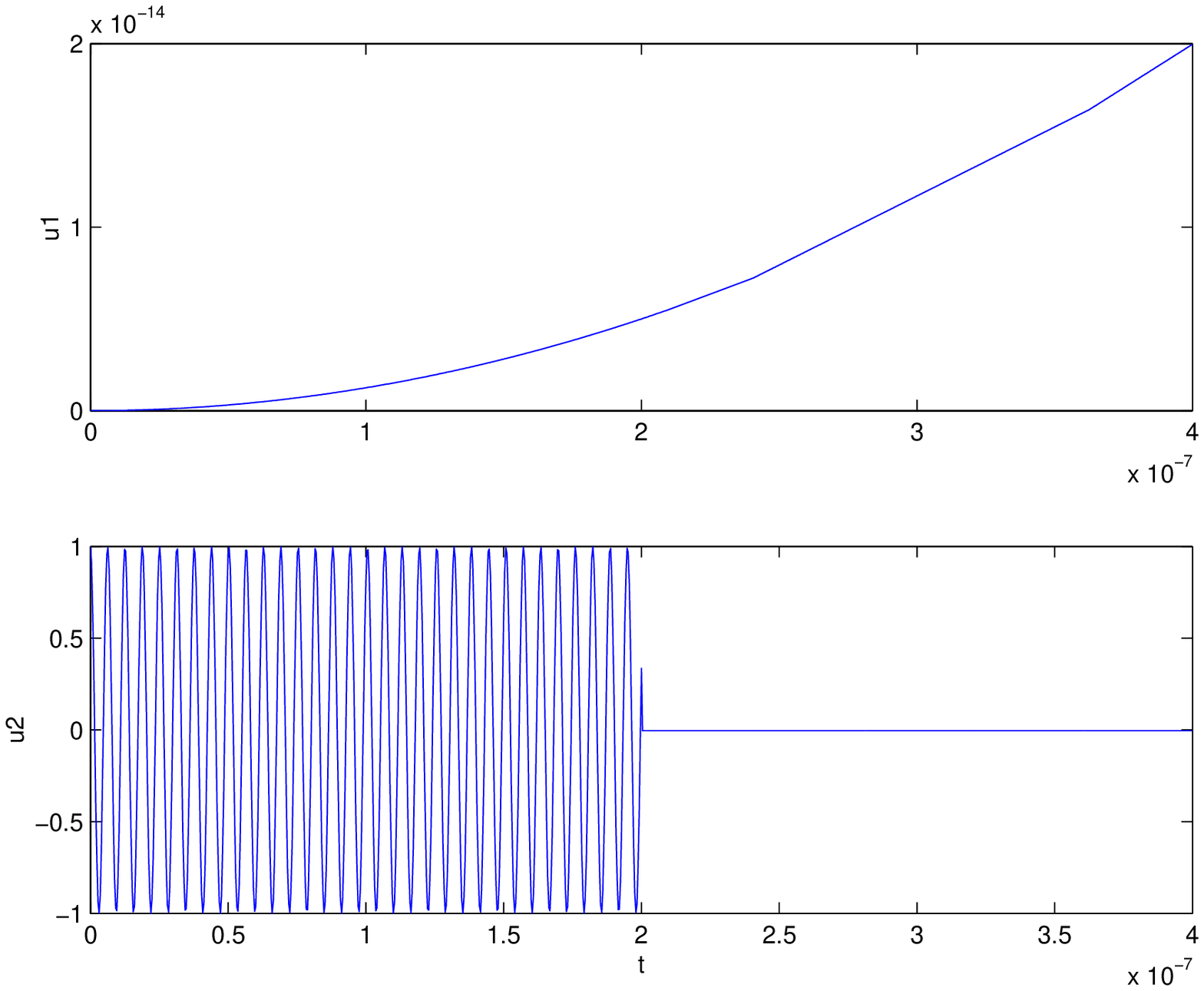}
    \caption{The solution of the simple model problem (\ref{eq:model})
      on $[0,100]$ (above) and on $[0,4\cdot{10}^{-7}]$ (below). The automatic
      modeling is activated at time $t = 2\tau = 2\cdot{10}^{-7}$.}
    \label{fig:model,solution}
  \end{center}
\end{figure}

\subsection{A lattice with internal vibrations}

\begin{figure}[htbp]
  \begin{center}
    \psfrag{m}{$m$}
    \psfrag{M}{$M$}
    \includegraphics[width=8cm]{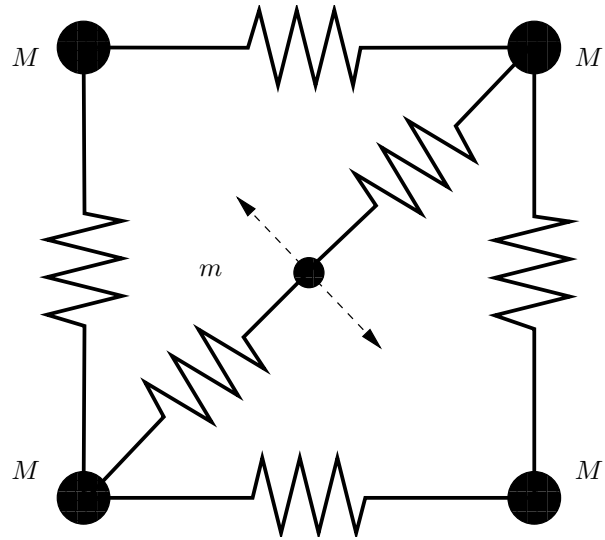}
    \caption{Detail of the lattice. The arrows indicate the direction
    of vibration perpendicular to the springs connecting the small
    mass to the large masses.}
    \label{fig:latticedetail}
  \end{center}
\end{figure}

\begin{figure}[htbp]
  \begin{center}
    \includegraphics[width=8cm]{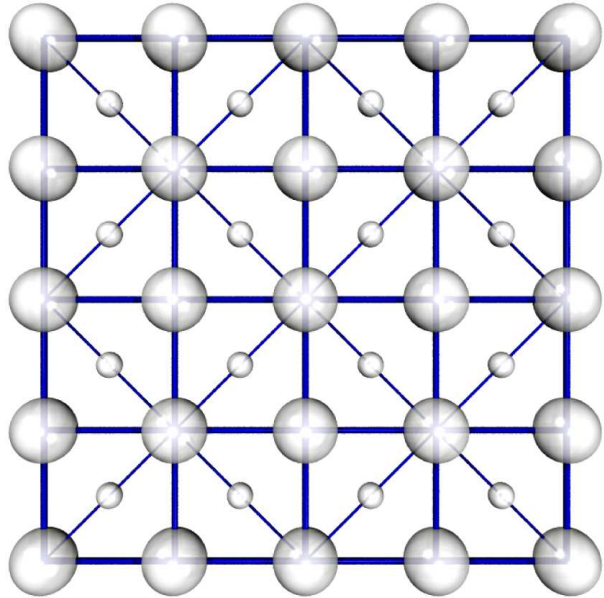}
    \caption{Lattice consisting of $p^2$ large masses and $(p-1)^2$ small masses.}
    \label{fig:lattice}
  \end{center}
\end{figure}

The second example is a lattice consisting of a set of $p^2$ large
and $(p-1)^2$ small point masses connected by springs of equal stiffness $\kappa = 1$,
as shown in Figure \ref{fig:latticedetail} and Figure \ref{fig:lattice}.
Each large point mass is of size $M = 100$ and each small point mass is
of size $m = 10^{-12}$, giving a large time scale of size $\sim 10$ and a
small time scale of size $\sim 10^{-6}$.

The fast oscillations of the small point masses make the initially
stationary structure of large point masses contract. Without resolving
the fast time scales and ignoring the subgrid model, the distance $D$
between the lower left large point mass at $x = (0,0)$ and the upper
right large point mass at $x = (1,1)$ remains constant, $D = \sqrt{2}$.
In Figure \ref{fig:lattice,solution2}, we show the computed solution
with $\tau = 10^{-4}$, which manages to correctly capture the oscillation
in the diameter $D$ of the lattice as a consequence of the internal
vibrations at time scale $10^{-6}$.

With a constant subgrid model $\bar{g}$ as in the example, the reduced
model stays accurate until the configuration of the lattice has
changed sufficiently. When the change becomes too large, the reduced
model can no longer give an accurate representation of the full
system, as shown in Figure \ref{fig:lattice,solution1}. At this point,
the reduced model needs to be reconstructed in a new short time
simulation.

\begin{figure}[htbp]
  \begin{center}
    \psfrag{t}{$t$}
    \psfrag{d}{$D(t)$}
    \includegraphics[width=10cm]{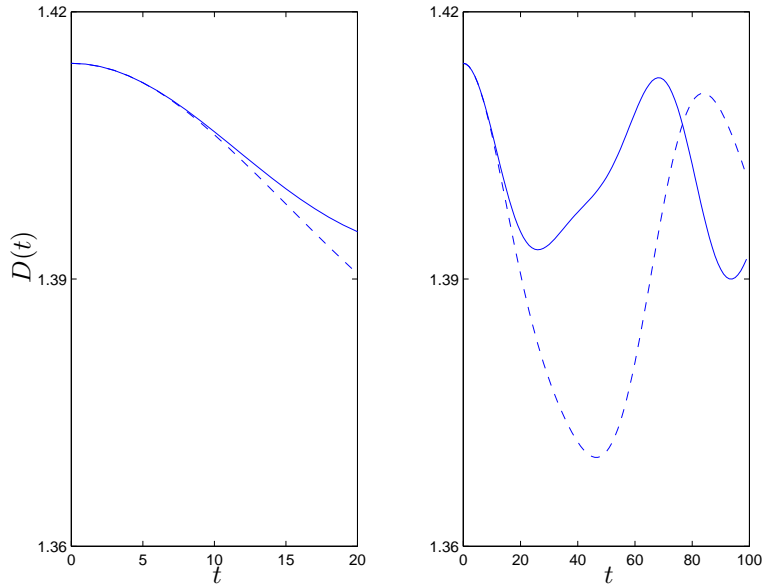}
    \caption{The diameter $D$ of the lattice as function of time on
      $[0,20]$ (left) and on $[0,100]$ (right) for $m = 10^{-4}$ and
      $\tau = 1$. The solid line represents the diameter for the solution of the
      reduced system (\ref{eq:reduced}) and the dashed line
      represents the solution of the full system (\ref{eq:u'=f}).}
    \label{fig:lattice,solution1}
  \end{center}
\end{figure}

\begin{figure}[htbp]
  \begin{center}
    \psfrag{t}{$t$}
    \psfrag{d}{$d(t)$}
    \psfrag{D}{$D(t)$}
    \includegraphics[width=10cm]{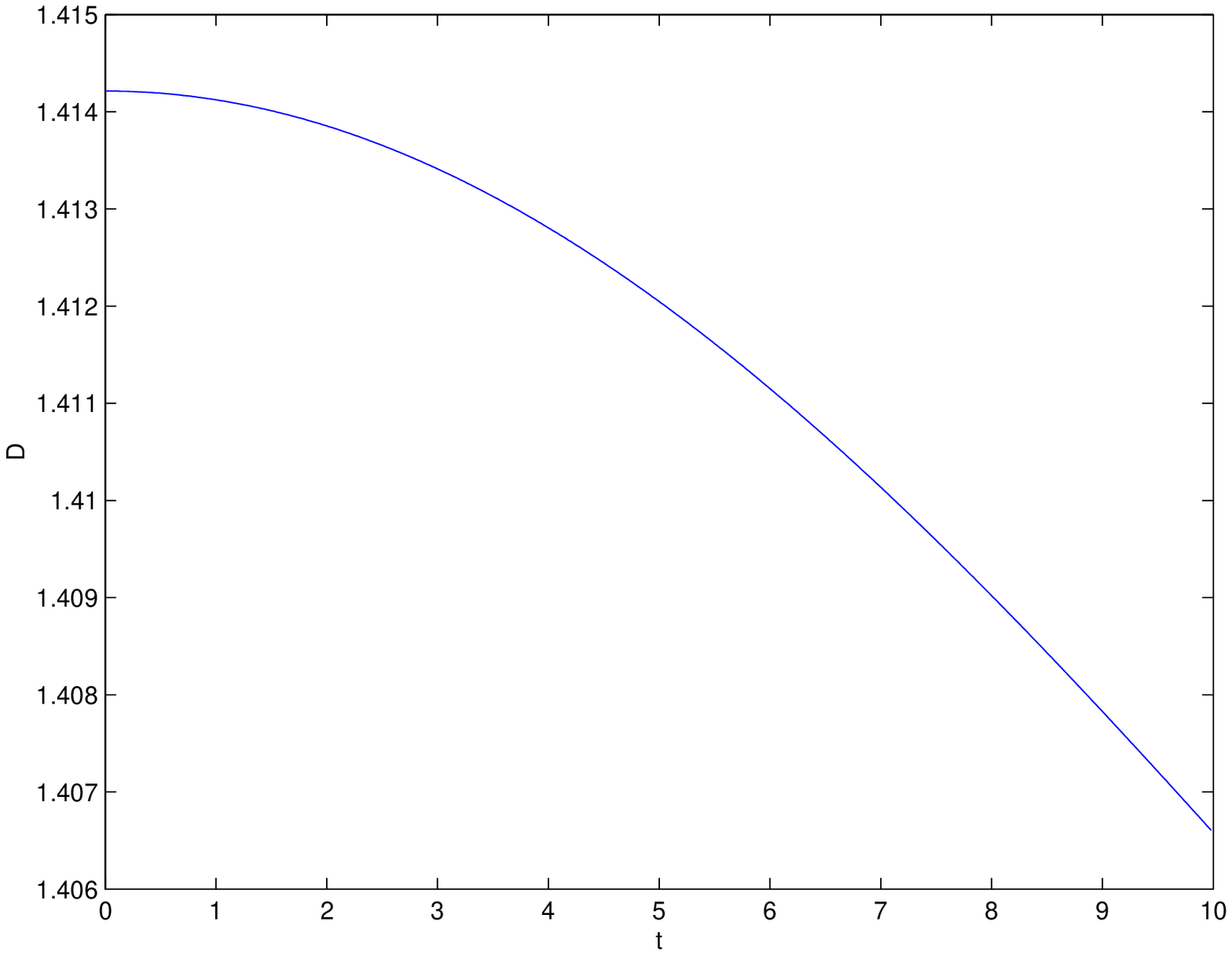} \\
    \includegraphics[width=10cm]{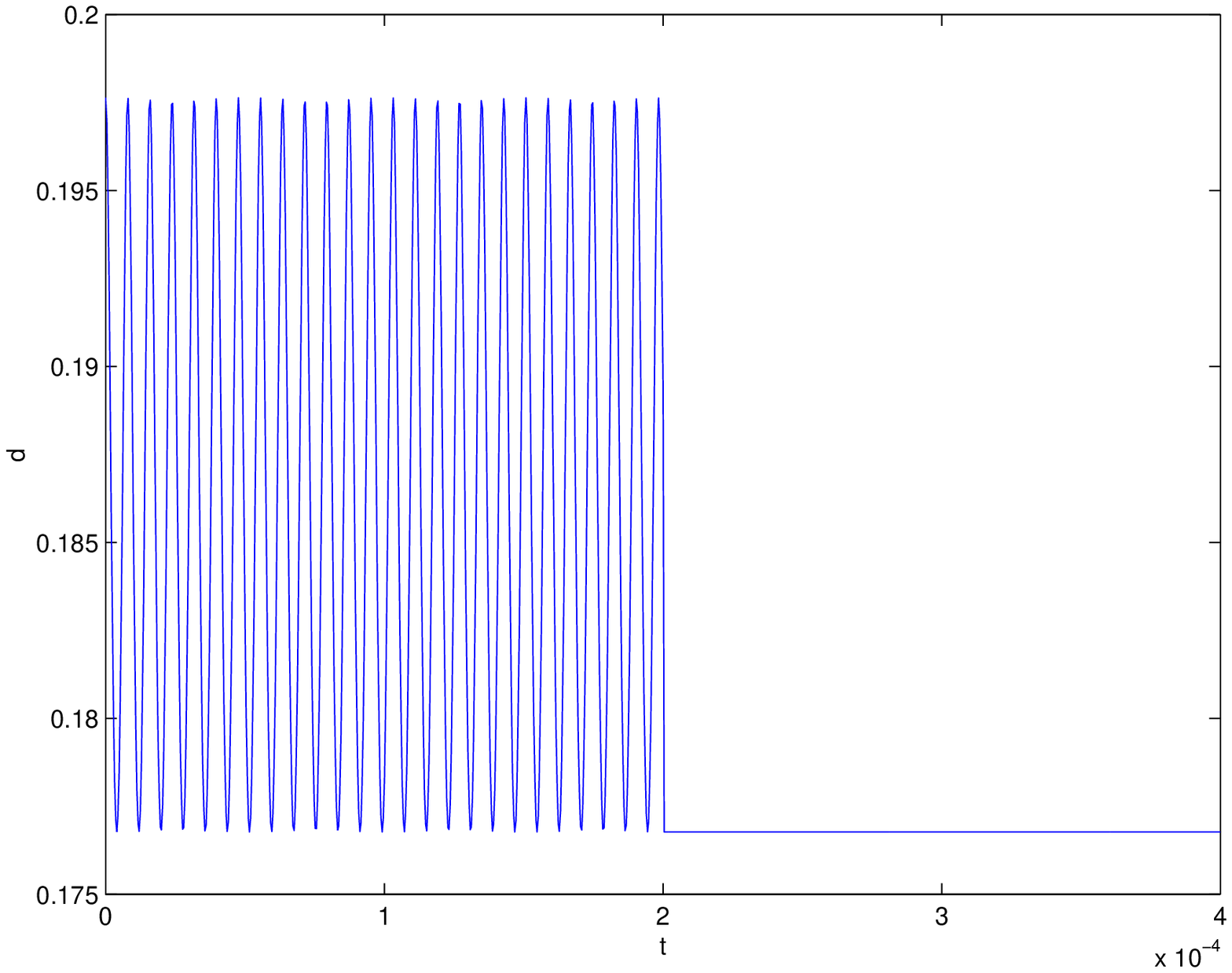}
    \caption{Distance $D$ between the lower left large mass and the upper
      right large mass (above) and the distance $d$ between the lower
      left large mass and the lower left small mass (below) as function of
      time on $[0,10]$ and on $[0,4\cdot 10^{-4}]$, respectively.}
    \label{fig:lattice,solution2}
  \end{center}
\end{figure}

\bibliographystyle{siam}
\bibliography{bibliography}

\end{document}